\documentclass{article}
\usepackage{amsmath,amssymb,amsthm}
\usepackage{mathrsfs}
\usepackage[colorlinks=true]{hyperref}
\usepackage{graphicx,subfigure}

\renewcommand{\geq}{\geqslant}
\renewcommand{\leq}{\leqslant}

\newtheorem{theorem}{Theorem}
\newtheorem{lemma}{Lemma}
\newtheorem{proposition}{Proposition}
\newtheorem{corollary}{Corollary}
\newtheorem{remark}{Remark}

\newcommand{\OCP}{\bf (OCP)}

\title{Optimal Neumann control for the 1D wave equation: Finite horizon, infinite horizon,
boundary  tracking terms and the turnpike property}

\author{Martin Gugat\footnote{FAU, Friedrich-Alexander-University Erlangen-N\"urnberg,
Department Mathematik, Cauerstr. 11, 91058 Erlangen, Germany (\texttt{martin.gugat@fau.de})}
\and
Emmanuel Tr\'elat\footnote{Sorbonne Universit\'es, UPMC Univ Paris 06, CNRS UMR 7598, Laboratoire Jacques-Louis Lions, Institut Universitaire de France, F-75005, Paris, France (\texttt{emmanuel.trelat@upmc.fr}).}
\and
Enrique Zuazua\footnote{BCAM - Basque Center for Applied Mathematics, Mazarredo, 14 E-48009 Bilbao-Basque Country-Spain.
Ikerbasque, Basque Foundation for Science, Alameda Urquijo 36-5, Plaza Bizkaia, 48011, Bilbao-Basque Country-Spain (\texttt{zuazua@bcamath.org})}
}

\date{}

\begin{document}

\maketitle

\begin{abstract}
We consider a vibrating string that is fixed at one end  with Neumann control action at the other end. We investigate the optimal control problem of steering this system from given initial data to rest, in time $T$, by minimizing an objective functional that is the convex sum of the $L^2$-norm of the control and of a boundary Neumann tracking term.

We provide an explicit solution of this optimal control problem, showing that if the weight of the tracking term is positive, then the optimal control action is concentrated at the beginning and at the end of the time interval, and in-between it decays exponentially.
We show that the optimal control can actually be written in that case as the sum of an exponentially decaying term and of an exponentially increasing term.
This implies that, if the time $T$ is large
%, then,
 %except at the beginning and the end of the time interval, the optimal control and the %corresponding state are exponentially close to $0$, and then
 the optimal trajectory approximately consists of three arcs, where the first and the third short-time arcs are transient arcs, and in the middle arc
 the  optimal control and the corresponding state are exponentially close to $0$.
 %is
 %in long time.
 This is   an example for  a turnpike phenomenon for a problem of optimal boundary control.
If $T=+\infty$ (infinite horizon time problem), then only the exponentially decaying component of the control remains, and the norms of the optimal control action and of the optimal state decay exponentially in time.
In contrast
 to this situation
 if the weight of the tracking term is zero
  and only the control cost is minimized, then the optimal control  is distributed uniformly along the whole interval $[0,T]$ and coincides with the
  %usual
  control given by the Hilbert Uniqueness Method.

%Moreover, in a similarity theorem, we establish that the solutions of the first problem and of the third one are similar, in the sense that, for every $T\in \{2,4,6,...\}$, there exists a weight $\lambda$ such that the controls coincide along the time interval $[0,2]$.
\end{abstract}


%\begin{AMS}
%35L05, 93C20, 35L20.
%\end{AMS}

\section{Introduction}\label{sec_introduction}
The turnpike property has been discussed recently for the
optimal control of linear systems governed by
partial differential equations, see \cite{porretta}.
Turnpike theory has originally been discussed in economics,
see \cite{samuelson}.
For systems governed by  ordinary differential equations the turnpike property
has been discussed for example in \cite{trelat}, also for the nonlinear case.
The turnpike property states loosely speaking
that if the
objective function penalizes
both control cost and the
difference of the optimal trajectory to
a given desired stationary state,
the optimal controls
will steer the system quickly to the
desired stationary state and
then the system will remain on this path most
of the time.
In Section 4 of \cite{porretta}, optimal control problems with the wave equation are considered
on a given finite time interval $[0,T]$.
The control is distributed in the interior of the domain
and no conditions for the terminal state at the time $T$ are prescribed.
In this paper we consider problems of optimal boundary control of
the wave equation.
We consider both problems with infinite time horizon and problems with finite time horizon.
For the finite time problems,
we prescribe exact terminal conditions.
If the objective function only penalizes the control cost,
for  the $1D$ case in \cite{gu:neumann} it has been shown that the
optimal controls are periodic, so in particular
they do not have the turnpike structure.
This illustrates that the turnpike property heavily depends on
the choice of the objective functions,
that has to couple the control cost with the
penalization of the distance to the desired state.

We consider a system governed by the one-dimensional wave equation on a finite space interval, with a homogeneous Dirichlet boundary condition at one side, and a Neumann boundary control action at the other:
\begin{equation}\label{waveeq}
\begin{split}
& \partial_{tt} y(t,x) = \partial_{xx} y(t,x), \qquad (t,x)\in \mathbb{R} \times (0,1), \\
& y(t,0)= 0,\ \partial_x y(t,1)= u(t), \qquad t\in \mathbb{R},
\end{split}
\end{equation}
where the control $u$ belongs to the class of square-integrable functions.

Let $y_0 \in H^1(0,1)$ be such that $y_0(0)=0$, and let $y_1 \in L^2(0,1)$ be arbitrary. For every $u\in L^2(0,+\infty)$, there exists a unique solution $y\in C^0(0,+\infty;H^1(0,1))\cap C^1(0,+\infty;L^2(0,1))$ of \eqref{waveeq} such that $y(0,\cdot)=y_0(\cdot)$ and $\partial_t y(0,\cdot)=y_1(\cdot)$.

As is well known, the system is exactly controllable if and only if $T\geq 2$ (see for example \cite{guleskly}). In this paper, given any $\lambda\in [0,1]$ and any $T\in [0,+\infty]$, we consider the optimal control problem $\OCP_\lambda^T$ of finding a control $u\in L^2(0,T)$ minimizing the objective functional
\begin{equation}\label{defJ}
J_\lambda^T(u) = \int_0^T \left( (1-\lambda) \left(\partial_x y(t,0) \right)^2 + \lambda u(t)^2 \right) dt ,
\end{equation}
such that the corresponding solution of \eqref{waveeq}, with $y(0,\cdot)=y_0(\cdot)$ and $\partial_t y(0,\cdot)=y_1(\cdot)$, satisfies $y(T,\cdot) = \partial_t y(T,\cdot) = 0$ at the final time (exact null controllability problem). If $T=+\infty$, then one can drop the final constraint requirement, which, by the way, happens to be automatically satisfied in the sense of a limit (this is a well known result coming from Riccati theory).

\medskip

For $\lambda=1$, this very classical optimal control problem (minimization of the $L^2$ norm of the control) has been considered, e.g., in \cite{Krabs,Lions,russell} (see also \cite{colombo} for optimal control problems consisting of satisfying consumer demands at the boundary of the system, such as in gas transportation networks), and it is easy to see that the optimal control is periodic (see \cite{gu:neumann}), with a period equal to $4$, that is, twice the time needed by a wave starting at the boundary point where the control acts to return to that point.
Note that, in this case, the optimal control is as well given by the Hilbert Uniqueness Method (see \cite{Lions}).
%At the uncontrolled boundary, the norm of the space derivative of the optimal state solution %may then be large if $T$ is large.

\medskip

If $\lambda<1$, then the objective functional involves a nontrivial boundary tracking term.
This tracking term may be considered as a boundary observation of the space derivative of the state at the uncontrolled end of the string.
As we are going to prove, in that case, the optimal control action is then essentially concentrated at the beginning and at the end of the time interval $[0,T]$. More precisely, the optimal control can be written as the sum of an exponentially decaying term and of an exponentially increasing term.

As a consequence, if $T$ is large then the optimal control, solution of $\OCP_\lambda^T$, approximately consists of three pieces: the first and the third pieces are in short-time, and are transient arcs; the middle arc is in long time, and is exponentially close to $0$. This is a turnpike phenomenon (see \cite{porretta,trelat}), meaning that the optimal trajectory, starting from given initial data, very quickly approximately reaches the steady-state $(0,0)$ (within exponentially short time, say $\varepsilon>0$), then remains exponentially close to that steady-state within long time (say, over the time interval $[\varepsilon,T-\varepsilon]$), and, in the last short-time part $[T-\varepsilon,T]$, leaves this neighborhood in order to quickly reach its target.

In this approximate picture, if $T=+\infty$ (infinite horizon), then the last transient arc does not exist since the infinite-horizon target is the steady-state $(0,0)$. In that case, the norm of the optimal control decays exponentially in time, and the same is true for the optimal state. Indeed, smallness of the observation term for a sufficiently long time interval with zero control implies proportional smallness of the state (this follows from an observability inequality, see \cite{zuazua:sirev}).

Another possible picture illustrating the turnpike behavior is the following. For $T$ large, the optimal trajectory of $\OCP_\lambda^T$ approximately consists of three arcs: the first arc is the solution of $\OCP_\lambda^\infty$ (infinite horizon problem), forward in time, and converges exponentially to $0$. The second arc, occupying the main (middle) part of the time interval, is the steady-state $0$. The third arc is the solution of $\OCP_\lambda^\infty$, but backward in time.
Note that the optimal control problem $\OCP_\lambda^\infty$
%in infinite horizon
fits into the well known Linear Quadratic Riccati theory.

\medskip

In all cases, we will provide completely explicit formulas for the optimal controls, which explain and imply the turnpike behavior observed for $\lambda<1$. This is in contrast with the case $\lambda=0$ for which the control action is distributed uniformly along the time interval $[0,T]$.
In addition, we will also establish a similarity theorem showing that, for every $T$ that is a positive even integer, there exists an appropriate weight $\lambda<1$ for which the optimal solutions of $\OCP_\lambda^T$ and of $\OCP_\lambda^\infty$ coincide along $[0,2]$.

%%%%%%%%%%%%%%%%%%%%%%%%%%%%%%%%%%%%%%%%%%%%%%
%%%%%%%%%%%%%%%%%%%%%%%%%%%%%%%%%%%%%%%%%%%%%%

In this paper we focus on problems governed by  the 1D wave equation.
In order to illustrate the generality of the turnpike phenomenon,
% illustrate that the turnpike structure is a general phenomenon,
% that is ubiquitous,
before we turn to the 1D wave equation
%that
% observed for the optimal controls for the 1D  wave equation
%
%also occurs in a more general framework,
we consider an example  in a more general framework with a  strongly continuous semigroup in
the following section.

%%%%%%%%%%%%%%%%%%%%%%%%%%%%%%%%%%%%%%%%%%%%%%%%%%%%%%%%%%%%%%%%%%%%%%%%%%%%%

\section{A general remark}

Let a Hilbert space $X$ be given.
Let $A$ be the generator of
a strongly continuous
semigroup $({\mathbb T}_t)_{t\geq 0}$
(for the definitions see
for example
\cite{lasiecka},
%Chapter 10 in
\cite{tu:obser}).
 Let $U$ be another Hilbert space that contains the controls.
 Let the linear operator  $B:  U\rightarrow  X$ be given.
Let a time $T>0$,  a weight $\gamma>0$ for the control cost
 and an initial state $y_0\in X$ be given.
 For $u\in L^2((0,T); U)$ and $t\in [0,T]$,
 define
 %\begin{equation}
 \[
 \Phi_t(u) =
 %(\beta I- A)
 \int_0^t {\mathbb T}_{t-\sigma}
 %(\beta I - A)^{-1}
 B u(\sigma)\, d \sigma.
\]
% \end{equation}
 We consider a system that is governed by the differential equation
 $y'(t)= A \, y(t) + B \, u(t)$,   $t\in [0,T]$ with the initial condition
 $y(0)  = y_0\in X$.
 We assume that there exists a time $T_{\min}>0$  such that for all
 $T\geq T_{\min}$
  the considered system is
 exactly controllable   in the sense that
 \[{\rm Ran}\,  \Phi_T=X.\]
 Let us assume that $T\geq T_{\min}$, then
 there exists a control function $u \in L^2((0,T); U)$ such
 that the terminal constrain $y(T)=0$ holds.
 Let $\lambda \in (0,1)$ be given.
 Consider the problem of optimal exact control
 \begin{equation}
 \label{08052015}
 \left\{
 \begin{array}{l}
  \min
 \int_0^T  (1 - \lambda) \langle y(t), y(t) \rangle_X +
  \lambda \langle u(t), u(t) \rangle_U\, dt
 %\\
 %\min \gamma \int_0^T \langle u(s), u(s) \rangle_U\, ds
 %+ \int_0^T \langle y(s), y(s) \rangle_X\, ds
 \\
  {\rm subject\;\; to}\;\;
 y'(t)= A \, y(t) + B \, u(t),\, t\in [0,T];
 \\
 y(0)  = y_0,\,  y(T) = 0.
 \end{array}
 \right.
 \end{equation}
%Let us point out that
The {\em static} optimal control problem corresponding to
(\ref{08052015}) is
 \begin{equation}
 \label{08052015a}
 \left\{
 \begin{array}{l}
 \min  (1 - \lambda) \langle y, y \rangle_X
 +\lambda  \langle u, u \rangle_U
 \\
  {\rm subject\;\; to}\;\;
 0 = A \, y + B \, u.
 \end{array}
 \right.
 \end{equation}
Obviously the solution of (\ref{08052015a}) is zero.
The solution of the static  optimal control problem (\ref{08052015a})
determines the turnpike which in our case is
$(u_{tp},y_{tp})=(0,0)$.
Results about the convergence of the long time average of
the optimal controls to the turnpike control can be found in
\cite{carlson}, \cite{zaslavski}.

 In order to determine the structure of the optimal control
 that solves (\ref{08052015})
 we look at the necessary optimality conditions.
For all $\delta_1 \in L^2((0,T); U)$,
 $ p \in C((0,T); X)$,
 $ \delta_2  \in C((0,T); X)$ with $\delta_2(0)=0$, $\delta_2(T)=0$
we have
\[
0 =  \int_0^T \lambda  \,\langle u(s),\,\delta_1(s) \rangle_U
+  (1 - \lambda) \langle y(s),\,\delta_2(s) \rangle_X
+
 \langle \delta_2'(s) - A\, \delta_2(s) - B \, \delta_1(s),\, p(s) \rangle_X
\, ds.
\]
This yields the optimality system
\begin{equation}
\left\{
\begin{array}{rrr}
%\begin{eqnarray*}
y' & = & A y + Bu,
\\
p' & = & - A^\ast p + y ,
\\
u  & = &
%\frac{1}{\gamma}
\frac{1 - \lambda}{\lambda}
 B^\ast \, p
\end{array}
\right.
\end{equation}
with the conditions $y(0)=y_0$, $y(T)=0$.
Hence we get
\begin{equation}
\label{hilf20150505}
- A p' = A  A^\ast p - A y = A  A^\ast p - y' +
%\frac{1}{\gamma}
\frac{1 - \lambda}{\lambda}
B  B^\ast \, p.
\end{equation}
By taking the time derivative in the first order equation for $p$ we get
the second order equation
%\[
$p''   =   - A^\ast p' + y'$.
%\]
This implies
\begin{equation}
\label{ystrichgleichung2015}
y' = p''   +  A^\ast p'.
\end{equation}
Now we can use
(\ref{ystrichgleichung2015}) to eliminate $y'$ from (\ref{hilf20150505})
and get
\[
- A p' =  A  A^\ast p
-p''   -  A^\ast p'
 +
 %\frac{1}{\gamma}
 \frac{1 - \lambda}{\lambda}
  B  B^\ast \, p.
 \]
 This yields
\begin{equation}
\label{DGL05052015}
p''   = (A  A^\ast  +
%\frac{1}{\gamma}
\frac{1 - \lambda}{\lambda}
 B  B^\ast)  \, p  + ( A - A^\ast) \, p' .
\end{equation}

%%%%%%%%%%%%%%%%%%%%%%%%%%%%%%%%%%%%%%%%%%%%%%%%%%%%%%%%%%%%%

\subsection{Skew-adjoint operators}
Now we consider the  case where $A$ is skew-adjoint, that is $A^\ast = -A$.
In \cite{porretta}, turnpike inequalities for the  case of the wave equation
where $A^\ast=-A$ are given in  Section 4.
%This corresponds to the case of  hyperbolic systems, that is
%discussed in \cite{porretta} in Section
Equation (\ref{DGL05052015}) yields
\begin{equation}
\label{DGL05052015b}
p''   =   L\, p + 2 A \, p'
\end{equation}
where
$ L = A  A^\ast  +
%\frac{1}{\gamma}
\frac{1 - \lambda}{\lambda}\,
 B  B^\ast$.

If there exist  solutions $D_+$, $D_-$ of the operator equation
\[D^2 = L + 2 A D\]
this yields  solutions of the form
\[p(t)= \exp(D_+\,t)\, p_1 + \exp(D_-\,t) p_2\]
where $p_1$, $p_2\in X$ are chosen such that
for the  state
$y  =   p' + A^\ast p$
we have
$y(0)=0$ and $y(T)=0$.
Note that $L+ A^2 =
\frac{1 - \lambda}{\lambda}
%\frac{1}{\gamma}
B  B^\ast $ is positive
in the sense that
$\langle x, (L+A^2) x \rangle_X \geq 0$ for all $x\in X$.
If
%$LA = AL$,
$A \, B\, B^\ast$ is  skew adjoint
(for example if $B$ is the identity)
%or $B\, B^\ast =0$),
we have $ D_{\pm}=  A \pm (L + A^2)^{1/2} $.

Now we assume that  $A$  and $B\,B^{\ast}$ are diagonalizable with the same sequence of orthonormal
eigenfunctions $(\varphi_k)_k$ and that the  real parts of the eigenvalues of $B\,B^{\ast}$ are bounded from below by $\omega^2>0$.
Then (\ref{DGL05052015b}) yields  a sequence of  ordinary differential equations
for $h_k(t)= \langle p(t), \varphi_k\rangle_X$ namely
\begin{equation}
h_k'' = \langle \varphi_k,\, L \,  \varphi_k\rangle_X \,  h_k +  2 \langle  \varphi_k,\, A \,  \varphi_k\rangle_X \, h_k'
.
%=: \lambda_k \,  h_k + 2  \mu_k  \, h_k'.
\end{equation}
With the  roots $\delta_k^+$, $\delta_k^-$ of the  characteristic polynomial
\[p_k(z)= z^2 - 2 \langle  \varphi_k,\, A \,  \varphi_k\rangle_X  \, z -  \langle \varphi_k,\, L \,  \varphi_k\rangle_X \]
we get the solutions
\[h_k(t) = u_k \exp(\delta_k^- t) + v_k  \exp(\delta_k^+ t)\]
where $ {\rm  Re}(\delta_k^+) \geq  \omega>0$ and $ {\rm  Re}(\delta_k^-) \leq  -\omega <0$.
%Note that $|{\rm  Re}(\delta_\pm)| > \omega>0.   $
The  coefficients $u_k$ and $v_k$ are chosen such that
$p'(0) = y_0 +  A p(0)$ and $ p'(T) = A p(T)$, because then we have $y(0)=y_0$ and $y(T)=0$.
In fact, this yields a constant $C_{\min}$ 
that is independent of $y_0$ and $T$
%only depends on $T_{\min}$ 
such that for all $T\geq T_{\min}$,
we have the inequality 
\[\sum_k \left(|u_k|^2 +   |\exp(2\,\delta_k^+ T)|\, |v_k|^2 \right)
\leq
%\frac{1}{ 1 - \exp(-2 \omega \, T)}
\,
C_{\min} \, \|y_0\|^2_X.
%  \, \frac{\|y_0\|^2}{\omega^2}.
\]
Using Parseval's equation we get the inequality
\begin{eqnarray}
\label{pturnpike18115}
\|p(t)\|_X
& \leq &
e^{-{\omega} t}
\;
\left(\sum_k |u_k|^2\right)^{1/2}
+
e^{-{\omega} (T -t)}
\left(\sum_k   |\exp(2\,\delta_k^+ T)|\, |v_k|^2)\right)^{1/2}
\\
& \leq &
\left(e^{-{\omega} t}  + e^{-{\omega} (T -t)}\right)\,
2\,\sqrt{C_{\min}} \|y_0\|_X
 % \, \frac{\|y_0\|}{\omega}
\label{pturnpike18115b}
%&
%+
%&
%+
%\left(
% \|\exp(D_+ \, T)  p_2\|
 %\right)
.
\end{eqnarray}
%\end{equation}
%Note that
%%for all $T\geq T_{\min}$,
%we have an inequality of the form
%\[\sum_k \left(|u_k|^2 +   \exp(2\,\delta_k^+ T) |v_k|^2 \right)
%\leq
%\frac{1}{ 1 - \exp(-2 \omega \, T)}
%  \, \frac{\|y_0\|^2}{\omega^2}.
%\]
% Korrekt aber zu lang.
%with a constant
%$C$  that is independent of $T$.

%
%Now let us assume that $D_+$ and $D_-$  are diagonalizable
%and the  real parts of the eigenvalues of $D_+$ are bounded from below by $\omega>0$
%and that
%the  real parts of the eigenvalues of $D_-= - D_+^\ast$ are bounded above by $-\omega<0$.
%%
%Then we have the inequality
%\begin{equation}
%\label{pturnpike18}
%\|p(t)\|
%\leq
%e^{-{\omega} t}
%\;
%\|p_2\|
%+
%e^{-{\omega} (T -t)}
%%\\
%%\label{pturnpike28}
%%&
%%+
%%&
%%+
%%\left(
% \|\exp(D_+ \, T)  p_2\|
% %\right)
%.
%%\end{eqnarray}
%\end{equation}
%%
%%where
%%$D$ solves
%%that is
%%$D = A - (L + A^2)^{1/2}$.
Inequality (\ref{pturnpike18115})-(\ref{pturnpike18115b}) is
a {\em turnpike inequality} for $p$.
It states that the norm of $p(t)$ is bounded above
%be decomposed as a sum of
by a sum of  a part that is
exponentially decreasing with time and
a second part that is exponentially increasing
towards $T$.
The  optimal control has the form
$u  =
%\frac{1}{\gamma}
\frac{1 - \lambda}{\lambda}
 B^\ast \, p$,
 so it also shows a turnpike structure.
  Note that the optimal control norms are
decreasing with $T$, hence they are uniformly bounded.
%So in this sense, also here $p$ shows
%a turnpike structure.

%This yields a solution of the form
%\begin{equation}
%\label{hyperbolicfunctionsa}
%p(t) = \cosh(t \, L^{\tfrac{1}{2}}) u_c +
%L^{- \tfrac{1}{2}} \, \sinh(t \, L^{\tfrac{1}{2}}) u_s
%+ 2\int_0^t L^{-\tfrac{1}{2}} \, \sinh( (t-s) L^{\tfrac{1}{2}}) A p'(s)\, ds
%\end{equation}
%with the $\cosh$ and $\sinh$ operators as defined in
%\cite{fattorini} and $u_c$, $u_s\in U$.
%

%%%%%%%%%%%%%%%%%%%%%%%%%%%%%%%%%%%%%%%%%%%%%%%%%%%%%%%%%%%%%%%%%%%%%%%%%%%%%%%%%%%%%%%%%%

\subsection{Self-adjoint operators}
Now we consider the case that $A$ is self-adjoint.
In \cite{porretta}, turnpike inequalities for the parabolic case
where $A^\ast=A$ are given in  Section 3.
Equation
(\ref{DGL05052015}) yields
\begin{equation}
\label{DGL05052015a}
p''   =   L\, p
\end{equation}
where
$ L = A  A^\ast  +  \frac{1}{\gamma} B  B^\ast$.
Thus we get
\[p(t)= \cosh(t\, L^{\tfrac{1}{2}}) \, p(0) + L^{- \tfrac{1}{2}} \, \sinh(t \, L^{\tfrac{1}{2}})\, p'(0)\]
with the $\cosh$ and $\sinh$ operators as defined in
\cite{fattorini}.
For the optimal state we have
\begin{eqnarray*}
y  & = &  p' + A^\ast p
%\\
%& = &
%  A^\ast \cosh(t\, L^{\tfrac{1}{2}})  p(0)
%  +\cosh(t\, L^{\tfrac{1}{2}}) p'(0)
%  \\
%  & + &
% A^\ast \, L^{- \tfrac{1}{2}}  \sinh(t \, L^{\tfrac{1}{2}})\, p'(0)
% +  \,L^{\tfrac{1}{2}} \sinh(t \, L^{\tfrac{1}{2}})\, p(0)
.
\end{eqnarray*}
The equations $y(0)=y_0$, $y(T)=0$ yield a system of linear equations for
$p(0)$, $p'(0)$.
In particular we have
$p'(0)=y_0 - A^{\ast} \, p(0)$.

Now let us assume that $L$ is diagonalizable
and the eigenvalues of $L$ are bounded from below by $\omega^2>0$.
%To be precise,
Then we have the turnpike inequality
\begin{eqnarray*}
\label{pturnpike1}
\|p(t)\|
&\leq &
 \frac{1}{2}
\exp(-{\omega} t)
\;
(\|p(0) \| + \|L^{-1/2} p'(0)\|)
+\frac{1}{2}
\exp(-{\omega} (T -t))
\\
\label{pturnpike2}
&
&
( \|\exp(L^{1/2} T) p(0)\|
+
 \|\exp(L^{1/2} T)  L^{-1/2} p'(0)\|)
.
\end{eqnarray*}
The  optimal control has the form
\begin{equation}
\label{hyperbolicfunctions}
u(t)   = \frac{1}{\gamma} B^\ast\cosh(t \, L^{\tfrac{1}{2}}) p(0) + \frac{1}{\gamma} B^\ast L^{- \tfrac{1}{2}} \, \sinh(t \, L^{\tfrac{1}{2}}) p'(0).
\end{equation}
This
%implies
means that also
%in this example
the optimal control
can be represented as the sum
of families of increasing and decreasing exponentials
with rates ${\omega} $.

\section{The main results}\label{sec_mainresults}
Now we come to our results about
optimal control problems for a  system governed by (\ref{waveeq}).
Let $y_0 \in H^1(0,1)$ be such that $y_0(0)=0$, and let $y_1 \in L^2(0,1)$ be arbitrary.

\subsection{Explicit formulas}
We have the following result, giving the explicit solution of $\OCP_\lambda^T$, for any $\lambda\in[0,1]$ and any $T\in[2,+\infty]$. Note that, when $T=+\infty$, we have to assume that $\lambda<1$ to ensure well-posedness.

\medskip

\begin{theorem}\label{thmpt}
For every $T\in [2,+\infty]$ and every $\lambda\in[0,1]$, the problem $\OCP_\lambda^T$ has a unique optimal control solution denoted by $u_\lambda^T$.
\begin{enumerate}
\item
We assume that $T=2n$ for some $n\in \mathbb{N}^*$.
\begin{itemize}
\item If $\lambda=1$, then the optimal control $u_1^T$, solution of $\OCP_1^T$ is $2$-anti-periodic, and thus $4$-periodic, meaning that
\begin{equation}\label{periodicneumann}
u_1^T(t+2k) = (-1)^k u_1^T(t),
\end{equation}
for every $t\in (0,2)$ and for every $k\in \mathbb{N}$ such that $t+2k\leq T$, and moreover,
\begin{equation}\label{exqactcontrolrepresentation}
u_1^T(t) = \left\{ \begin{array}{ll}
\frac{1}{2n}  \left(y_0'(1-t) - y_1(1-t)\right), &  t\in (0,1), \\[2mm]
\frac{1}{2n}  \left(y_0'(t-1) + y_1(t-1)\right), & t\in (1,2).
\end{array} \right.
\end{equation}

\item If $\lambda<1$, then the optimal control solution of $\OCP_\lambda^T$ is the sum of an  exponentially decaying term and of an exponentially increasing one. More precisely, defining the real number $z_\lambda\in (-1,0)$ by
\begin{equation}\label{zgammadefinition}
z_\lambda = \frac{-\lambda}{2-\lambda + 2\sqrt{1-\lambda}} ,
\end{equation}
we have
\begin{equation}\label{ugammatrepresentation}
u_\lambda^T(t+2k) =  z_\lambda^k f_+(t) + \frac{1}{z_\lambda^k} f_-(t) ,
\end{equation}
for every $t\in (0,2)$ and every $k\in \mathbb{N}$ such that $t+2k\leq T$, where
\begin{equation}\label{f+-def}
f_+(t) = \frac{1 + z_\lambda}{ 1 - z_\lambda^{2n}}\, F(t-1),\qquad
f_-(t) = \frac{1 + \frac{1}{z_\lambda} }{ 1 - \frac{1}{z_\lambda^{2n}}}\, F(t-1),
\end{equation}
with $F\in L^2(-1,1)$ defined by
\begin{equation}
\label{Fdefinition}
F(t) = \left\{ \begin{array}{ll}
\frac{1}{2} \left( y_0'( - t) -  y_1(-t) \right), & t \in (-1,0),\\ [2mm]
\frac{1}{2} \left( y_0'(t)  +  y_1(t) \right), & t \in [0,1).
\end{array} \right.
\end{equation}
\end{itemize}

\item We assume that $T=+\infty$ and that $\lambda<1$.

If $\lambda=0$, then the optimal control $u_0^\infty$, solution of $\OCP_0^\infty$, coincides along the time interval $[0,2]$ with the optimal control $u_1^2$, solution of $\OCP_1^2$.

If $0<\lambda<1$, the optimal control $u_\lambda^\infty$, solution of $\OCP_\lambda^\infty$, is given along the time interval $[0,2]$ by
\begin{equation}\label{ugammadefinition}
u_\lambda^\infty(t) = \left\{ \begin{array}{ll}
\frac{  1 + z_\lambda}{2}  \left( y_0'(1-t)  - y_1(1-t)\right), &  t\in (0,1), \\ [2mm]
\frac{ 1 + z_\lambda}{2}  \left( y_0'(t-1)  + y_1(t-1)\right) , & t\in (1,2),
\end{array}
\right.
\end{equation}
and moreover, we have
\begin{equation}\label{uzkgammagleichung}
u_\lambda^\infty(t +  2 k) =  z_\lambda^k  u_\lambda^\infty(t),
\end{equation}
for every $t\in (0,2)$ and every $k\in \mathbb{N}^*$.

The corresponding optimal state $y_\lambda^\infty$ decays exponentially, in the sense that there exists $C_0>0$ such that
\begin{multline}\label{decayinfinity}
\int_0^1 \left( \left( \partial_x y_\lambda^\infty(t + 2\, k,\, x) \right)^2 +  \left( \partial_t y_\lambda^\infty(t + 2 \, k,\, x) \right)^2 \right) dx \\
\leq C_0 |z_\lambda|^{2k} \int_0^1 \left(  y_0'(x)^2 + y_1(x)^2\right) dx ,
\end{multline}
for every $t\in (0,2)$ and every $k\in \mathbb{N}^*$.
\end{enumerate}
\end{theorem}

\medskip

Theorem \ref{thmpt} is proved in Section \ref{proof_thmpt}.

\medskip

For $\lambda=1$, that is, when there is no tracking term in the objective functional, the explicit solution of $\OCP_1^T$ given above has already been computed in \cite[Theorem 2.1]{gu:neumann}.
In this case, the problem consists of minimizing the $L^2$ norm of the (Neumann) control. The optimal control $u_1^T$, whose explicit formula is given above, can also be characterized as well with the famous Hilbert Uniqueness Method (see \cite{Lions}) and is then often referred to as the \emph{HUM control}.

Here, there is no dissipation induced by the objective functional (no tracking term), the optimal control is periodic, and is uniformly distributed over the time interval $[0,T]$, in the sense that there is no energy decay.

\medskip

In contrast, if $\lambda<1$, the control is the sum of two terms, one of which is exponentially decreasing, and the other being exponentially increasing. For $T$ large enough, this implies the turnpike phenomenon, stated in details in Section \ref{sec_turnpike}.

\medskip

\begin{remark}
For $\lambda=0$, the solution of $\OCP_0^\infty$ coincides with the solution of the problem of optimal feedback control studied in \cite{gu:optimal}.
\end{remark}

\medskip

\begin{remark}
The estimate (\ref{decayinfinity}) is clearly equivalent to
$$
\int_0^1 \left( \left( \partial_x y_\lambda^\infty(t,\, x) \right)^2 +  \left( \partial_t y_\lambda^\infty(t,\, x) \right)^2 \right) dx \leq C_1 e^{- \mu  t} \int_0^1 \left( y_0'(x)^2 + y_1(x)^2 \right) dx ,
$$
for every time $t\geq 0$, for some positive constants $C_1$ and $\mu$ not depending on the initial state (see also \cite[Lemma 2]{gutus}).
\end{remark}

\medskip

\begin{remark}
It is well known that the solution of the infinite horizon problem $\OCP_\lambda^\infty$ can also be expressed in feedback form (Linear Quadratic Riccati theory),
see for example \cite{flandoli}. More precisely, the velocity feedback
$$
\partial_x y(t,1)= \frac{z_\lambda+1}{z_\lambda -1}\, \partial_t y(t,1)
$$
generates the same state as the one generated by the optimal control $u_\lambda^\infty$.
\end{remark}

\medskip

\begin{remark}\label{rem_steadystate}
In the above results, we considered only the null steady-state, but we can easily replace it with any other steady-state, as follows. Any steady-state of \eqref{waveeq} is given by $\bar y(x) = \sigma x$, for some $\sigma \in \mathbb{R}$. Then, all results therein can be written in terms of such a steady-state: it suffices to replace, everywhere, $y(t,x)$ with $y(t,x)-\sigma x$, and $\partial_x y(t,x)$ with $\partial_x y(t,x)-\sigma$. For instance, the right boundary condition becomes $\partial_x y(t,1)=\sigma+u(t)$, the final conditions become $y(T,x)=\sigma x$ and $\partial_t y(T,x)=0$, and the objective functional becomes
$$
J_\lambda^T(u) = \int_0^T \left( (1-\lambda) \left(\partial_x y(t,0)-\sigma \right)^2 + \lambda u(t)^2 \right) dt .
$$
\end{remark}

\subsection{Consequence: the turnpike behavior}\label{sec_turnpike}
From Theorem \ref{thmpt} and from the previous discussions, we infer the following consequence on the qualitative behavior of the optimal solution.

\medskip

\begin{corollary}\label{thm_turnpike}
For every $\lambda\in[0,1)$, then there exist $C_1>0$ and $\mu>0$ such that, for every $T\geq 2$, for all initial conditions $(y_0,y_1)\in H^1(0,1)\times L^2(0,1)$ with $y_0(0)=0$, the optimal solution of $\OCP_\lambda^T$ satisfies the estimate
\begin{equation}\label{estim_turnpike}
\int_0^1 \left( \left( \partial_x y_\lambda^T(t,\, x) \right)^2 +  \left( \partial_t y_\lambda^T(t,\, x) \right)^2 \right) dx \leq C_1 e^{- \mu t(T-t)} \int_0^1 \left( y_0'(x)^2 + y_1(x)^2 \right) dx ,
\end{equation}
for every $t\in[0,T]$.
\end{corollary}

\medskip

In the estimate \eqref{estim_turnpike}, what is important to see is that the term $e^{- \mu t(T-t)}$ is equal to $1$ at times $t=0$ and $t=T$, but it is exponentially small in the middle of the interval. It becomes even smaller and smaller when $T$ is taken larger. This estimate implies the turnpike behavior described previously: short-time arcs at the beginning and at the end of the interval are devoted to satisfy the terminal constraints, and in-between, the trajectory remains essentially close to rest.

%The turnpike property has been investigated  in \cite{porretta}
%for problems of optimal distributed control and in the general finite-dimensional nonlinear %case in \cite{trelat} (see also references therein).

\subsection{Similarity result}
We next state the following similarity result: for any final time $T$ that is a positive even integer, there exists a weight $\lambda$ such that the optimal solutions of $\OCP_1^\infty$ and $\OCP_\lambda^\infty$ coincide along the subinterval $[0,2]$ of $[0,T]$.

\medskip

\begin{theorem}\label{similaritytheorem}
Given any $T\in 2 \mathbb{N}^*$, we choose $\lambda>0$ such that
\begin{equation}\label{gammachoice}
z_\lambda = \frac{2}{T} - 1.
\end{equation}
Then we have
\begin{equation}\label{identischesteuerung}
u_1^T(t) =u_\lambda^\infty(t),\quad \forall t\in (0,2),
\end{equation}
and
\begin{equation}\label{similarityungleichung}
\left\| u_1^T(\cdot) - u_\lambda^\infty(\cdot) \right\|_{L^2(2k,2k+2)} \leq  \left| 1 - |z_\lambda|^k \right| \, \frac{2}{T} \left( \| y_0' \|_{L^2(0,1)} + \| y_1\|_{L^2(0,1)}\right).
\end{equation}
for every $k\in \{0,1,\ldots,(T-2)/2\}$.
\end{theorem}

\medskip

The proof of Theorem \ref{similaritytheorem} is done in Section \ref{proof_similaritytheorem}.

\subsection{Numerical illustration}
We set $y_0(x) = 4\sin(\pi x/2)$ and $y_1(x)=0$, for every $x\in[0,1]$.
From Theorem \ref{thmpt}, if $0<\lambda<1$ then the optimal control solution of $\OCP_\lambda^T$ is given by
$$
u_\lambda^T(t  + 2 k ) = |z_\lambda|^k \, (1 + z_\lambda)\, \pi \,\sin\left(\frac{\pi}{2}( t + 2 k) \right),
$$
for $t\in (0,2)$ and $k\in \mathbb{N}$.

\begin{figure}[h]
\begin{center}
\subfigure[$\lambda=24/25$ and $T=20$.]
{\includegraphics[width=6.25cm]{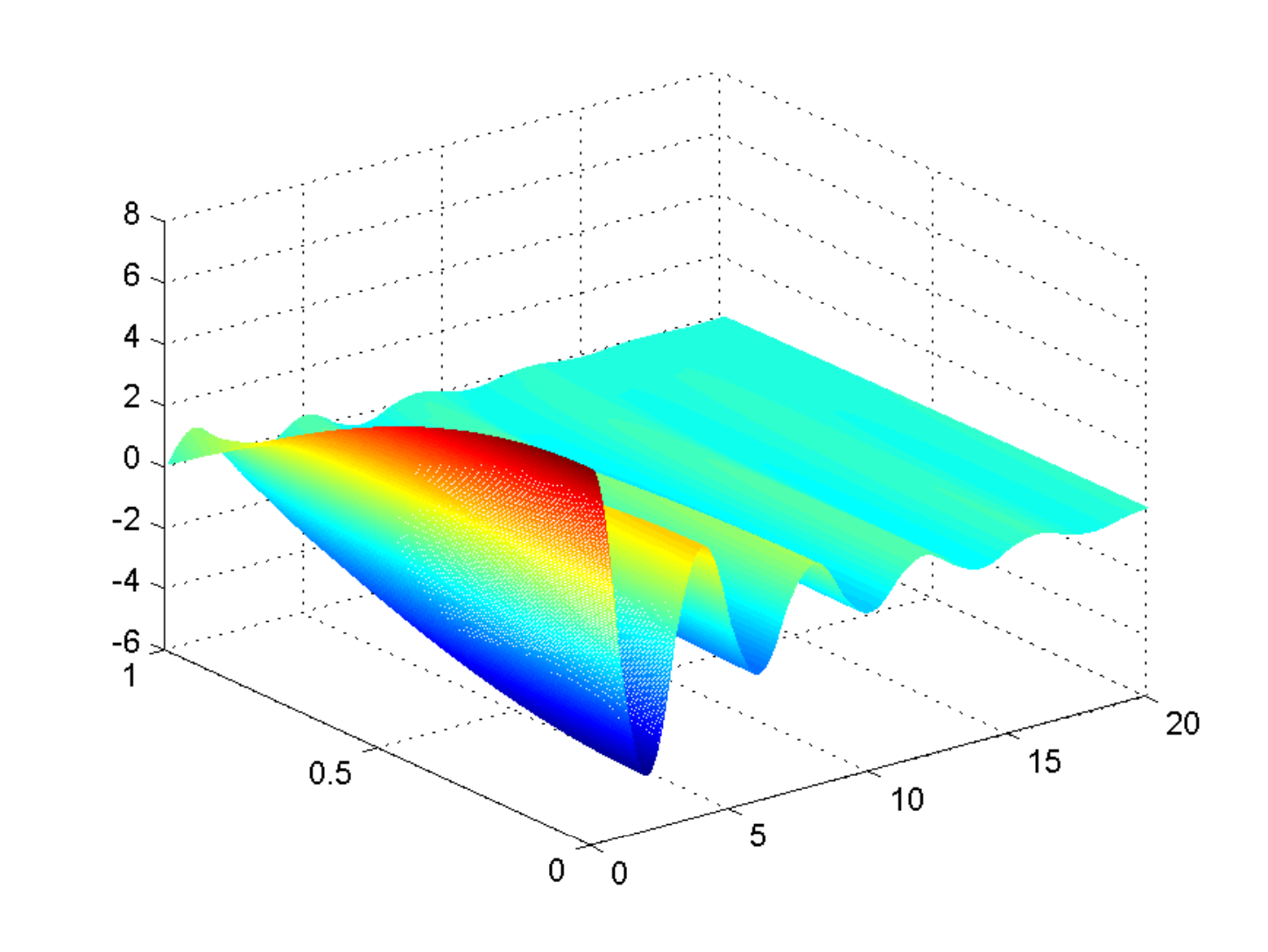}\label{bild1}}
%\qquad\qquad\qquad
\subfigure[$\lambda=99/100$ and $T=20$.]
{\includegraphics[width=6.25cm,height=4.7cm]{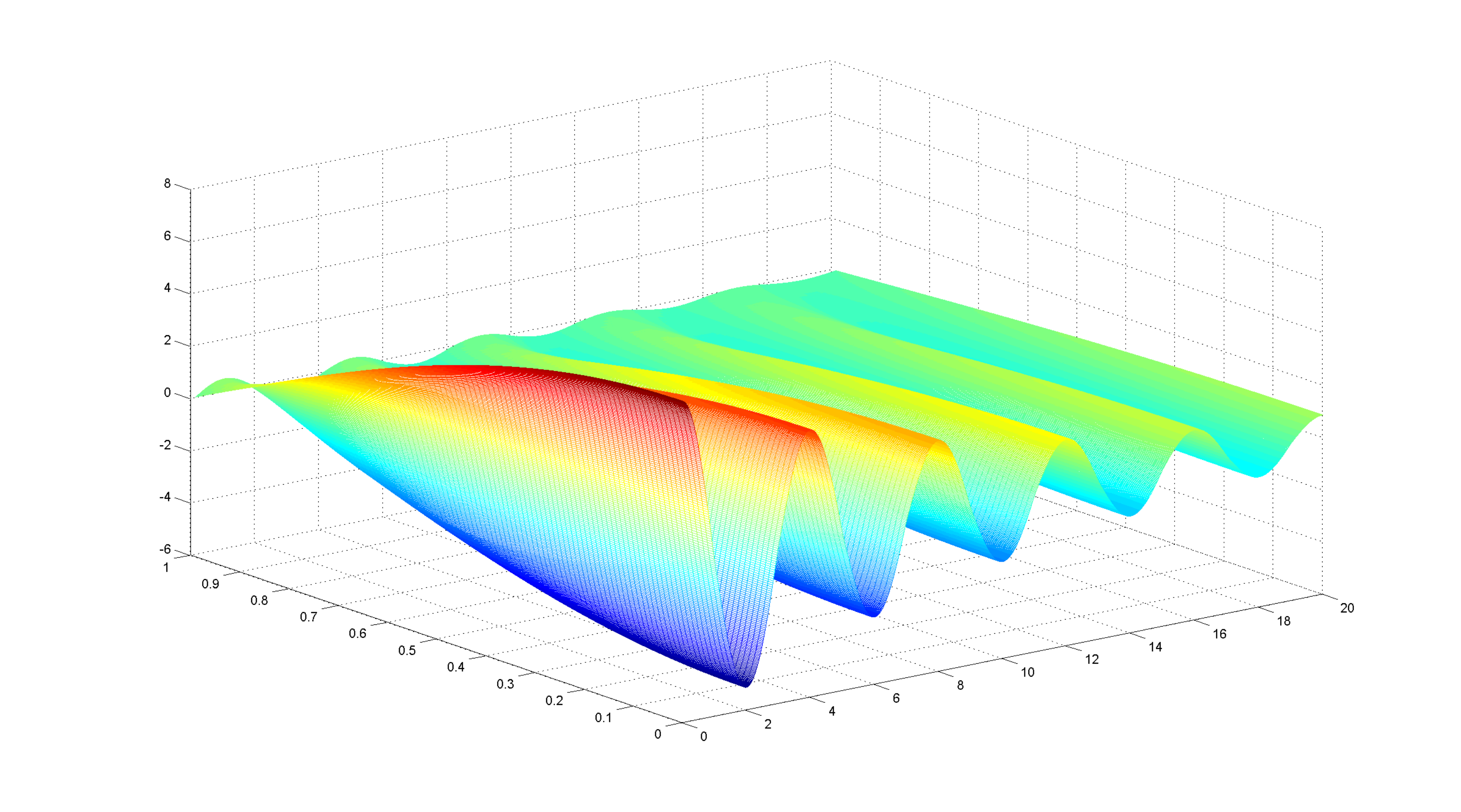}\label{bild2}}
\caption{Plot of $\partial_x y(t,x)$.} \label{bild}
\end{center}
\end{figure}

The graph of $\partial_x y_\lambda^T(t,x)$ is provided on Figure \ref{bild}, for $T=20$, with $\lambda=24/25$ on Figure \ref{bild1} and $\lambda=99/100$ on Figure \ref{bild2}.
The control $u_\lambda^T(t)=\partial_x y_\lambda^T(t,1)$ is the boundary trace at the back.

%\begin{figure}[h]
%\centerline{\includegraphics[width=7cm]{bildyxT=20gamma=24}}
%\caption{Plot of $\partial_x y(t,x)$, with $\gamma=24$ and $T=20$.}\label{bild1}
%\end{figure}
%
%\begin{figure}[h]
%\centerline{\includegraphics[width=12cm]{bildyxT=20gamma=99}}
%\caption{
%The figure shows $\partial_x y(t,x)$ for the optimal state for $\gamma=99$, the control $u_\lambda^\infty(t)=\partial_x y(t,1)$ is the trace at the back. The time interval is $[0,\,T]=[0,20]$.}\label{bild2}
%\end{figure}

These figures illustrate that the norm of the optimal state decays faster if $\lambda$ is smaller, as expected. However, smaller values of $\lambda$ cause larger oscillations. Note that $z_\lambda = - 2/3$ if $\lambda = 24/25$, and $z_\lambda = - 9/11$ if $\lambda = 99/100$.
Moreover, as pointed out in \cite{gu:neumann}, if $T\in 2\mathbb{N}^*$ then $u_1^T(t+2k) = \frac{2}{T} \pi \sin\left(\frac{\pi}{2} (t + 2 k)\right)$, for all $t$ and $k$ such that $t+2k\in [0,T]$ (see also \cite[Figure 4]{gu:neumann} for the corresponding optimal state, with $T=10$, up to the factor $2\pi$).

%%%%%%%%%%%%%%%%%%%%%%%%%%%%%%%%%%%%%%%%%%%%%%%%%%%%%%%%%%%%%%%%%%%%%%%%%%%%%

\section{Proofs}\label{sec_proofs}

\subsection{Well-posedness of the  initial-boundary value problem}
Let $y_0 \in H^1(0,1)$ be such that $y_0(0)=0$, and let $y_1 \in L^2(0,1)$ be arbitrary. Let $T\geq 2$, and let $u\in L^2(0,T)$ be fixed.
As a preliminary result, we study the well-posedness of the initial boundary value problem \eqref{waveeq} for a fixed control $u$, and with the fixed initial data $(y_0,y_1)$. The analysis is similar to the one done in \cite{gutus}.

We search a solution given as the sum of two traveling waves, i.e.,
$$
y(t,x)= \alpha(x+t) + \beta(x-t),
$$
where the functions $\alpha$ and $\beta$ are to be determined from the initial data and from the boundary data. First of all, to match the initial conditions, we must have
\begin{eqnarray}\label{ab1}
\alpha(t) & = & \frac{1}{2} \left(y_0(t) + \int_0^t y_1(s)\, ds \right)+C_0, \\ \label{ab2}
\beta(t)  & = & \frac{1}{2} \left(y_0(t) - \int_0^t y_1(s)\, ds \right)-C_0 ,
\end{eqnarray}
where $C_0$ is a real number. Besides, the boundary condition $y(t,0)=0$ implies that
\begin{equation}\label{betagleichung}
\beta(-s) = - \alpha(s),
\end{equation}
for almost every $s>0$. The boundary condition at $x=1$ leads to $\alpha'( 1 + t ) = u(t) - \beta'(1-t)$, and integrating in time, we get
$$
\alpha( t + 1) =  \beta(1-t)  +  \int_0^t u(s)\, ds + \alpha(1) - \beta(1) .
$$
Using (\ref{ab1}) and (\ref{ab2}), we have $\alpha(1) - \beta(1) =\int_0^1 y_1(s)\, ds + 2 C_0$, and therefore, choosing $C_0 =  -\frac{1}{2} \int_0^1 y_1(s)\, ds$, we get $\alpha(1) - \beta(1)=0$ and
\begin{equation}\label{alphagleichung0}
\alpha( t + 1) =  \beta(1-t) + \int_0^t  u(s)\, ds .
\end{equation}
Using (\ref{betagleichung}), the values of $\alpha$ for $t\in (0,1)$, given by (\ref{ab1}), determine those of $\beta$ for $t\in(-1,0)$.
The values of $\beta$ for $t\in (0,1)$ are given by (\ref{ab2}).
Now, knowing $\beta$ on the interval $(-1,1)$, we deduce from (\ref{alphagleichung0}) the values of $\alpha$ on the interval $(1,3)$.

Using (\ref{betagleichung}), we get $\alpha(t+1)= - \alpha(t-1) + \int_0^t u(s)\, ds$, for $t\geq 1$, or equivalently,
\begin{equation}\label{alphagleichung1}
\alpha(t+2)= - \alpha(t) + \int_0^{t+1}  u(s)\, ds,
\end{equation}
for $t\geq 0$.

Using (\ref{alphagleichung1}) enables us to determine $\alpha$ iteratively: starting with $\alpha$ on the interval $(1,3)$, the values of $u(t)$ yield those of $\alpha$ on $(3,5)$, and then using (\ref{alphagleichung1}), we determine $\alpha$ on $(7,9)$, etc.

In order to express everything in terms of $\alpha$ only (without using $\beta$), we extend the domain of $\alpha$ so that it contains $(-1,0)$. We get the values of $\alpha$ on $(-1,0)$ by using (\ref{betagleichung}) for $s\in (0,1)$, which yields $\alpha(t)=-\beta(-t)$ for $t\in (-1,0)$ with the values of $\beta$ on $(0,1)$ given by (\ref{ab2}).
 Then, using (\ref{alphagleichung0}), we get $\alpha(t+2)= - \alpha(t) + \int_0^{t+1} u(s)\, ds$ for $t\in (-1,0)$. We have the following lemma.

\medskip

\begin{lemma}\label{lemmaalpha}
Let $y_0 \in H^1(0,1)$ be such that $y_0(0)=0$, and let $y_1 \in L^2(0,1)$ be arbitrary.
We set
\begin{equation}\label{c0defi}
C_0 =  -\frac{1}{2} \int_0^1 y_1(s)\, ds,
\end{equation}
and we define $\alpha\in L^2(-1,1)$ by
\begin{equation}\label{abminus}
\alpha(t)  = \left\{\begin{array}{lcl}
\frac{1}{2} \left(- y_0( - t) + \int_0^{-t} y_1(s)\, ds \right) + C_0 & \textrm{if} & t \in (-1,0),   \\[2mm]
\frac{1}{2} \left(y_0(t) + \int_0^t y_1(s)\, ds \right)+C_0 & \textrm{if} & t \in [0,1).
\end{array}\right.
\end{equation}
Let $T=2K\in \mathbb{N}^*$, and let $u\in L^2(0,T)$ be fixed.

The function $\alpha$, defined by iteration according to
\begin{equation}
\label{recursion}
\alpha(t + 2k)= - \alpha(t + 2 (k-1)) + \int_0^{t+ 2k -1} u(s)\, ds,
\end{equation}
for every $t\in (-1,1)$ and every $k\in \mathbb{N}$ such that $t< T + 1 - 2k$, is well defined on the interval $(-1,T+1)$, and belongs to $H^1(-1,T+1)$.
\end{lemma}

\medskip

{\bf Proof:}
From the construction, it is clear that $\alpha|_{(k-1,k)}\in H^1(k,k+1)$, for every $k\in\mathbb{N}$. To prove that $\alpha \in H^1(-1,T+1)$, it suffices to prove that $\alpha$ is continuous. Since $\alpha(0^+) = \alpha(0^-)=  C_0$, $\alpha$ is continuous at $t=0$.
Using (\ref{abminus}), $\alpha$ is continuous as well on $(-1,1)$.

At $t=1$, using (\ref{c0defi}) we get $\alpha(1^+)=\alpha(1^-)=\frac{1}{2} y_0(1)$, and hence $\alpha$ is continuous at $t=1$.
Then, at this step, we have obtained that $\alpha$ is continuous on $(-1,3)$.

We then proceed by induction. Let $k\in \mathbb{N}^*$. We assume that $\alpha$ is continuous on the interval $(-1, 1+2k)$. Then $\alpha((-1+2k)^-)= \alpha((-1+2k)^+)$.
Using (\ref{recursion}), we have
\begin{multline*}
\alpha((1 + 2k)-) = - \alpha((1 + 2 (k-1))^- ) + \int_0^{2k} u(s)\, ds
= - \alpha((-1 +  2k )^- ) + \int_0^{2k} u(s)\, ds \\
= - \alpha((-1 +  2k )^+ ) + \int_0^{2k} u(s)\, ds
= \alpha((-1 + 2 (k+1))^+)
= \alpha((1 + 2k)^+).
\end{multline*}
Since $\alpha$ is defined by (\ref{abminus}), we infer that $\alpha$ is continuous on $(-1,1 + 2(k+1))$ for $k+1\leq K$. Lemma \ref{lemmaalpha} is proved.
$\Box$

\medskip

Using Lemma \ref{lemmaalpha}, we are now in a position to compute the solution of the initial boundary value problem under consideration in this subsection.

\medskip

\begin{proposition}\label{satz1}
Let $y_0 \in H^1(0,1)$ be such that $y_0(0)=0$, and let $y_1 \in L^2(0,1)$ be arbitrary.
Let $T\in2\mathbb{N}^*$, and let $u\in L^2(0,T)$ be fixed.
We consider the function $\alpha$ defined in Lemma \ref{lemmaalpha} by (\ref{abminus}).
Then the solution of \eqref{waveeq}, associated with the control $u$ and with the initial data $(y_0,y_1)$, is given by
\begin{equation}
\label{dalembertloesung}
 y(t,x)= \alpha(  t + x ) - \alpha( t - x ),
\end{equation}
for all $(t,x) \in (0,T)\times (0,1)$.
\end{proposition}

\medskip

{\bf Proof:}
The construction of $\alpha$ implies that $y$, defined by (\ref{dalembertloesung}), is a solution of the initial boundary value problem under consideration. We conclude by Cauchy uniqueness.
$\Box$

%%%%%%%%%%%%%%%%%%%%%%%%%%%%%%%%%%%%%%%%%%%%%%%

\subsection{Proof of Theorem \ref{thmpt}}\label{proof_thmpt}
%In this section we prove Theorem \ref{thmpt}.

\subsubsection{Case $T<+\infty$}
Let $y_\lambda^T$ be the state generated by the control $u_\lambda^T$ defined in the theorem. Let us first prove that $y_\lambda^T$ satisfies the terminal constraints
\begin{equation}\label{gammaterminalconstraints}
y_\lambda^T(T,\cdot)= 0, \quad \partial_t y_\lambda^T(T,\cdot)=0.
\end{equation}
It suffices to prove that $\alpha'(z)= 0$, for $z \in \,(T-1,\, T+1)$.
From (\ref{dalembertloesung}), we have $y_\lambda^T(t,x)= \alpha_\lambda(  t + x ) - \alpha_\lambda( t - x )$, with $\alpha_\lambda$ defined by (\ref{abminus}). The definition (\ref{Fdefinition}) implies that $F(t)=\alpha_\lambda'(t) $. Hence, we have
\begin{multline*}
\alpha_\lambda'(t-1) =   \frac{1 - z_\lambda^{2n}}{1 - z_\lambda^{2n}}\, F(t-1)
= \left( \frac{1}{1 - z_\lambda^{2n}} + \frac{z_\lambda^{2n}}{z_\lambda^{2n}-1 } \right) F(t-1) \\
= \left( \frac{1}{1 - z_\lambda^{2n}} + \frac{1}{1 - z_\lambda^{-T} } \right) F(t-1)
= \frac{1}{1 + z_\lambda} \, f_+(t) + \frac{1}{1 + \frac{1}{z_\lambda}} \, f_-(t)  ,
\end{multline*}
where the last equality follows from (\ref{f+-def}).

By (\ref{recursion}) we have
\begin{equation}\label{rekursionableitungt}
\alpha_\lambda'(t +  1) = - \alpha_\lambda'(t -1) + u_\lambda^\infty (t),
\end{equation}
for $ t \in (0,2)$. Using (\ref{ugammatrepresentation}), this yields
$\alpha_\lambda'(t +  1) = - \alpha_\lambda'(t -1) + f_+(t) + f_-(t)$,
for $ t \in (0,2)$, and then, using (\ref{f+-def}),
\begin{equation*}
%\begin{split}
\alpha_\lambda'(t +  1) %& = - \frac{1}{1 + z_\lambda} \, f_+(t) - \frac{1}{1 + \frac{1}{z_\lambda}} \, f_-(t) + f_+(t) + f_-(t) \\
%&
= \frac{z_\lambda}{1 + z_\lambda} \,  f_+(t) + \frac{ \frac{1}{z_\lambda}} {1 + \frac{1}{z_\lambda}} \, f_-(t).
%\end{split}
\end{equation*}
By induction, thanks to (\ref{rekursionableitungt}) and (\ref{ugammatrepresentation}), this implies that
\begin{equation}\label{decayrekursiontt}
\alpha_\lambda'(t -  1 + 2k  ) =  \frac{z_\lambda^k}{1 + z_\lambda} \,  f_+(t) + \frac{\frac{1}{z_\lambda^k}}{1 + \frac{1}{z_\lambda}} \, f_-(t) ,
\end{equation}
for every $t\in (0,2)$ and every $k\in \mathbb{N}$ such that $2k \leq T$.
Taking $k=T/2$, we get
\begin{equation}\label{decayrekursionttN}
\alpha_\lambda'(t -  1 + T  ) =  \frac{z_\lambda^{n}}{1 + z_\lambda} \,  f_+(t) + \frac{\frac{1}{z_\lambda^{n}}}{1 + \frac{1}{z_\lambda}} \, f_-(t) .
\end{equation}
Using (\ref{f+-def}), we infer that
\begin{equation}\label{decayrekursionttN1}
\alpha_\lambda'(t -  1 + T  ) = \left(  \frac{z_\lambda^{n}}{1 - z_\lambda^{2n}} + \frac{\frac{1}{z_\lambda^{n}}}{1 - \frac{1}{z_\lambda^{2n}}} \right) F(t-1) =  0 ,
\end{equation}
and hence the state $y_\lambda^T$ satisfies the terminal conditions (\ref{gammaterminalconstraints}).

For a control of the form $u = u_\lambda^T + h$, the generated state is $y=y_\lambda^T + y_h$, where $y_h$ is the state generated by the perturbation control $h$, with the boundary conditions $y_h(t,0)=0$, $\partial_x y_h(t,1)=h(t)$, and null initial conditions.
We only consider variations $h$ for which $y_h(T,\cdot)=0$ and $\partial_t y_h(T,\cdot)=0$.
Using (\ref{dalembertloesung}), we have
$$
y_h(t,x)= \alpha_h(  t + x ) - \alpha_h( t - x ).
$$
Since $y_h(0,\cdot)=\partial_t y_h(0,\cdot)=0$, we must have $\alpha_h'(t-1)=\alpha_h(t)=0$.  Moreover, owing to the terminal constraints, we must have $\alpha_h'=0$ along $(T-1,T+1)$.

The value of the objective functional of $\OCP_\lambda^T$ is
\begin{multline*}
\int_0^T \left( (1-\lambda) (\partial_x y(t,0) )^2 + \lambda u(t)^2\right) dt \\
= \int_0^T \Big( (1-\lambda)\left( (\partial_x y_\lambda^T(t,0) )^2 + (\partial_x y_h(t,0) )^2 + 2 \partial_x y_\lambda^T(t,0) \partial_x y_h(t,0) \right) \\
+ \lambda \left( u_\lambda^T(t)^2 + h(t)^2  + 2 u_\lambda^T(t) h(t) \right) \Big) \, dt.
\end{multline*}
We consider the linear part
$$
L_\lambda(h)= 2\int_0^T \left(  (1-\lambda) \partial_x y_\lambda^T(t,0) \partial_x y_h(t,0)
+ \lambda  u_\lambda^T(t) h(t) \right) dt.
$$
Since $\partial_x y_\lambda^T(t,0)= 2 \alpha_\lambda'(t)$, $\partial_x y_h(t,0)= 2\alpha_h'(t)$, and
\begin{equation*}
\begin{split}
u_\lambda^T(t) & = \partial_x y_\lambda^T(t,1) = \alpha_\lambda'(  t + 1 ) + \alpha_\lambda'( t - 1)  , \\
h(t) & = \partial_x y_h(t,1) = \alpha_h'(  t + 1 ) + \alpha_h'( t - 1) ,
\end{split}
\end{equation*}
we get
\begin{equation*}
\begin{split}
L_\lambda(h)
& = \int_0^T \Big( 8 (1-\lambda) \alpha_\lambda'(t) \alpha_h'( t)  \\
& \qquad\qquad + 2\lambda \left( \alpha_\lambda'( t+ 1) + \alpha_\lambda'( t- 1) \right) \left( \alpha_h'( t+ 1) + \alpha_h'( t- 1)\right) \Big) \, dt \\
& = \sum_{j=0}^{T-1} \int_0^1 \Big( 8 (1-\lambda) \alpha_\lambda'( t + j) \alpha_h'( t+ j) \\
&  \qquad + 2\lambda \left( \alpha_\lambda'( t+ 1 + j) + \alpha_\lambda'( t- 1 + j) \right) \left( \alpha_h'( t+ 1 + j) + \alpha_h'( t- 1 + j)\right) \Big) \, dt \\
& = 8 (1-\lambda) \sum_{j=1}^{T-1} \int_0^1 \alpha_h'( t + j) \alpha_\lambda'( t + j) \, dt \\
& \quad  + 2 \lambda \sum_{j=1}^{T-2} \int_0^1 \alpha_h'( t + j) \left( \alpha_\lambda'( t + 2  + j) +  \alpha_\lambda'( t+ j)  \right) dt  \\
& \quad +  2 \lambda  \sum_{j=1}^{T-2} \int_0^1 \alpha_h'( t + j) \left(  \alpha_\lambda'( t+ j)  + \alpha_\lambda'( t - 2  + j) \right) dt \\
& = 2 \sum_{j=1}^{T-2} \int_0^1 \alpha_h'( t + j) \Big( 4 (1-\lambda) \alpha_\lambda'( t+ j)  \\
& \qquad\qquad\qquad\qquad\qquad\qquad \lambda \left( \alpha_\lambda'( t + 2  + j) +  \alpha_\lambda'( t - 2  + j) + 2 \alpha_\lambda'( t+ j) \right)  \Big) \, dt
\end{split}
\end{equation*}
Defining the characteristic polynomial by
\begin{equation}\label{pgammadefinition}
p_\lambda(z)=  \lambda z^2 + (4-2\lambda) z + \lambda ,
\end{equation}
we have $p_\lambda(z_\lambda)=0$ and $p_\lambda(1/z_\lambda)=0$. Using (\ref{decayrekursiontt}), we have
$$
\lambda \alpha_\lambda'( t + 2  + j) + (4 - 2 \lambda) \alpha_\lambda'( t+ j) + \lambda \alpha_\lambda'( t - 2  + j) = 0,
$$
for every $t\in (0,1)$ and every $j\in \{1,2,\ldots,T-2\}$.
This implies that $L_\lambda(h)= 0$.
Now, concerning the value of the objective functional of $\OCP_\lambda^T$, for any $h$ such that $y_h(T,\cdot)=0$ and $\partial_t y_h(T,\cdot)=0$, we infer that
$$
\int_0^{+\infty} \left( (1-\lambda) \left(\partial_x y(t,0) \right)^2 + \lambda u(t)^2\right) dt \geq \int_0^{+\infty} \left( (1-\lambda) \left(\partial_x y_\lambda^T(t,0) \right)^2 + \lambda u_\lambda^T(t)^2 \right) dt,
$$
with a strict inequality whenever $h\neq 0$.
It follows that $u_\lambda^T$ is the unique optimal solution of $\OCP_\lambda^T$, as soon as $\lambda>0$.
If $\lambda=0$, then the result also follows from the representation of $L_0(h)$. However, in this case the characteristic polynomial $p_0(z)= 4z$ has only one root given by $z_0=0$.
Theorem \ref{thmpt} is proved for $T<+\infty$.

\subsubsection{Case $T=+\infty$}
We are going to use the previously established well-posedness results.

Let $y_\lambda^\infty$ be the state generated by the control $u_\lambda^\infty$ defined in the theorem. For a control of the form $u = u_\lambda^\infty + h$, the generated state is $y=y_\lambda^\infty + y_h$, where $y_h$ is the state generated by the control $h$, with null initial conditions and with the boundary conditions $y_h(t,0)=0$ and $\partial_x y_h(t,1)=h(t)$. The value of the objective functional of $\OCP_\lambda^\infty$ is
\begin{multline*}
\int_0^{+\infty} \left( (1-\lambda) (\partial_x y(t,0) )^2 + \lambda u(t)^2\right) dt \\
= \int_0^{+\infty} \Big( (1-\lambda) \left( (\partial_x y_\lambda^T(t,0) )^2 + (\partial_x y_h(t,0) )^2 + 2 \partial_x y_\lambda^\infty(t,0) \partial_x y_h(t,0) \right)  \\
+ \lambda \left( u_\lambda^\infty(t)^2 + h(t)^2  + 2 u_\lambda^\infty(t) h(t) \right) \Big) \, dt .
\end{multline*}
We consider the linear part
$$
L_\lambda(h)= 2\int_0^{+\infty} \left(  (1-\lambda) \partial_x y_\lambda^\infty(t,0)  \partial_x y_h(t,0) + \lambda u_\lambda^\infty(t) h(t) \right) dt .
$$
Using (\ref{dalembertloesung}), we have $y_\lambda^\infty(t,x)= \alpha_\lambda(  t + x ) - \alpha_\lambda( t - x )$ and $y_h(t,x)= \alpha_h(  t + x ) - \alpha_h( t - x )$, with $\alpha_\lambda$ given by (\ref{abminus}), and $\alpha_h=0$ on $(-1,1)$.
It follows that $\partial_x y_\lambda^\infty(t,0)= 2 \alpha_\lambda'(  t )$, $\partial_x y_h(t,0)= 2 \alpha_h'(  t )$, and
\begin{equation*}
\begin{split}
u_\lambda^\infty(t) & = \partial_x y_\lambda^\infty(t,1) =  \alpha_\lambda'(  t + 1 ) + \alpha_\lambda'( t - 1)   ,\\
h(t) & = \partial_x y_h(t,1) = \alpha_h'(  t + 1 ) + \alpha_h'( t - 1),
\end{split}
\end{equation*}
and therefore,
\begin{equation*}
\begin{split}
L_\lambda(h)
& = \int_0^{+\infty} \Big( 8 (1-\lambda) \alpha_\lambda'( t) \alpha_h'( t)  \\
& \qquad\qquad + 2\lambda  \left( \alpha_\lambda'( t+ 1) + \alpha_\lambda'( t- 1) \right)  \left(\alpha_h'( t+ 1) + \alpha_h'( t- 1)\right) \Big) \, dt  \\
& = \sum_{j=0}^{+\infty} \int_0^1 \Big( 8 (1-\lambda) \alpha_\lambda'( t + j) \alpha_h'( t+ j) \\
& \quad + 2\lambda \left( \alpha_\lambda'( t+ 1 + j) + \alpha_\lambda'( t- 1 + j) \right) \left( \alpha_h'( t+ 1 + j) + \alpha_h'( t- 1 + j)\right) \Big) \, dt \\
& = 2 \int_0^1 \Big( 4 (1-\lambda) \alpha_h'( t ) \alpha_\lambda'( t )  \\
& \qquad\quad +   \lambda\left( \alpha_h'( t - 1 ) \left( \alpha_\lambda'( t+ 1 ) + \alpha_\lambda'( t- 1 ) \right) + \alpha_h'( t ) \left( \alpha_\lambda'( t+ 2) + \alpha_\lambda'( t) \right) \right) \Big) \, dt \\
& \quad + 2 \sum_{j=1}^{+\infty} \int_0^1 \alpha_h'( t + j) \Big(  4(1-\lambda) \alpha_\lambda'( t+ j)  \\
&\qquad\qquad\qquad\qquad\quad\quad\ + \lambda\left( 2 \alpha_\lambda'( t+ j) + \alpha_\lambda'( t - 2  + j)  +  \alpha_\lambda'( t + 2  + j) \right)  \Big) \, dt.
\end{split}
\end{equation*}
By Lemma \ref{lemmaalpha}, for $t\in (0,1)$ the values of $\alpha_h'(t-1)$ and $\alpha_h'(t)$ are determined from the initial data, and since they are equal to zero, we have $\alpha_h'(t-1)=\alpha_h(t)=0$. This yields
\begin{multline}\label{lgammarepresentation}
L_\lambda(h) =  2 \sum_{k=1}^{+\infty} \int_0^1 \alpha_h'( t + k)  \Big( 4 (1-\lambda) \alpha_\lambda'( t+ k) \\
+ \lambda\left( 2 \alpha_\lambda'( t+ k) + \alpha_\lambda'( t - 2  + k) + \alpha_\lambda'( t + 2  + k) \right) \Big)\, dt .
\end{multline}
If $\lambda>0$ then the roots of the characteristic polynomial $p_\lambda$ defined by (\ref{pgammadefinition}) are $z_\lambda$ and $\frac{1}{z_\lambda}$.
In particular, we have $p_\lambda(z_\lambda)=0$. Note that, by Lemma \ref{lemmaalpha}, for $t\in (0,1)$ the values of $\alpha_\lambda'(t-1)$ and of $\alpha_\lambda'(t)$ are determined from the initial data. By (\ref{recursion}), we have
\begin{equation}\label{rekursionableitung}
\alpha_\lambda'(t +  1) = - \alpha_\lambda'(t -1) + u_\lambda^\infty (t),
\end{equation}
for $ t \in (0,2)$. Using the representation (\ref{ugammadefinition}) of $u_\lambda^\infty(t)$ for $ t \in (0,1)$, and using (\ref{abminus}), we infer that $\alpha_\lambda'(t +  1) = z_\lambda  \alpha_\lambda'(t -  1)$ for $ t \in (0,1)$.
Similarly, using (\ref{recursion}), we have $\alpha_\lambda'(t +  2) = - \alpha_\lambda'(t ) + u_\lambda^\infty (t+1)$ for $ t \in (0,1)$.
Using the representation (\ref{ugammadefinition}) of $u_\lambda^\infty(t)$ for $ t \in (1,2)$, and using (\ref{abminus}), we infer that $\alpha_\lambda'(t +  2) = z_\lambda \alpha_\lambda'(t )$ for $t\in (0,1)$.
It follows that $\alpha_\lambda'(t +  1) = z_\lambda \alpha_\lambda'(t -  1)$ for $t\in (0,2)$.
By induction, using (\ref{rekursionableitung}), (\ref{uzkgammagleichung}) and (\ref{ugammadefinition}), we get that
\begin{equation}\label{decayrekursion}
\alpha_\lambda'(t -  1 + 2k  ) = z_\lambda^k \alpha_\lambda'(t -  1  ) ,
\end{equation}
for every $t\in (0,2)$ and every $k\in \mathbb{N}$.
Therefore, we have obtained that
$$
\lambda \alpha_\lambda'( t + 2  + k) + (4 - 2 \lambda)  \alpha_\lambda'( t+ k) + \lambda \alpha_\lambda'( t - 2  + k)
= p_\lambda(z_\lambda) \alpha_\lambda'( t - 2  + k) = 0 ,
$$
for every $t\in (0,1)$ and every $k\in \mathbb{N}^*$.
We conclude that $L_\lambda(h)= 0$. Concerning the value of the objective functional of $\OCP_\lambda^\infty$, we infer that
\begin{multline*}
\int_0^{+\infty} \left( (1-\lambda) \left(\partial_x y(t,0) \right)^2 + \lambda u(t)^2\right) dt \\
\geq \int_0^{+\infty} \left( (1-\lambda) \left(\partial_x y_\lambda^\infty(t,0) \right)^2 + \lambda u_\lambda^\infty(t)^2 \right) dt ,
\end{multline*}
with a strict inequality whenever $h\neq 0$. It follows that $u_\lambda^\infty$ is the unique optimal solution of $\OCP_\lambda^\infty$ for $\lambda>0$.
For $\lambda=0$ the result also follows from the representation (\ref{lgammarepresentation})
of $L_0(h)$, with the difference that, in this case, the characteristic polynomial is $p_0(z)= z$ having the unique root $z_0=0$.

The inequality (\ref{decayinfinity}) follows from (\ref{decayrekursion}), since for the optimal state we have the energy
\begin{multline*}
\int_0^1 \left( \partial_x y_\lambda^\infty(t + 2 k, x)  \right)^2 +  \left( \partial_t y_\lambda^\infty(t + 2 k, x) \right)^2 dx = \int_{t+2k-1}^{t+2k+1} \alpha_\lambda'(s)^2 \, ds \\
= \int_{t-1}^{t+1} \alpha_\lambda'(s + 2k)^2 \, ds
= |z_\lambda|^{2k} \int_{t-1}^{t+1} \alpha_\lambda'(s)^2 \, ds ,
\end{multline*}
for every $t\in (0,2)$ and every $k\in \mathbb{N}$. Theorem \ref{thmpt} is proved for $T=+\infty$.

\medskip

\begin{remark}
The computation of the solution with the characteristic polynomial $p_\lambda$ is related to techniques used for linear difference equations, or for finite-dimensional linear systems with tridiagonal matrices (see \cite{meurant}).
\end{remark}

\subsection{Proof of Theorem \ref{similaritytheorem}}\label{proof_similaritytheorem}
We assume that $T = 2n$. Using (\ref{gammachoice}), we have $\frac{1}{2n} = \frac{ 1 + z_\lambda}{2}$, and then, using (\ref{exqactcontrolrepresentation}) and the representation
(\ref{ugammadefinition}) of $u_\lambda^\infty$, we infer (\ref{identischesteuerung}).
We have
$$
\left( \int_0^2 u_\lambda^\infty(t)^2 \,dt \right)^{1/2}  \leq  (1 + z_\lambda) \left( \| y_0' \|_{L^2(0,1)} + \| y_1\|_{L^2(0,1)}\right).
$$
The inequality (\ref{similarityungleichung}) follows similarly, using (\ref{periodicneumann}) and (\ref{uzkgammagleichung}), since
$$
\left\| u_\lambda^T -u_\lambda^\infty\right\|_{L^2(2k,2k+2)} = \left| 1 - |z_\lambda|^k \right| \left(\int_0^2 u_\lambda^\infty(t)^2 \,dt \right)^{1/2}.
$$
Theorem \ref{similaritytheorem} is proved.

%%%%%%%%%%%%%%%%%%%%%%%%%%%%%%%%%%%%%%%%%%%%%%%%%%%%%%%%%%%%%%%%%%%%%%%%%%%%%

\section{Conclusion}
We have discussed the influence of the objective function and of the time horizon on optimal Neumann boundary controls for the 1D wave equation. If the objective function is the control norm ($\lambda=1$) and if the terminal state is prescribed exactly, then the control action is distributed uniformly over the whole time horizon, and coincides with the control given by the Hilbert Uniqueness Method, which is periodic.

In contrast, if the objective function involves an additional tracking term ($\lambda<1$), then the optimal control action is essentially concentrated at the starting time $0$ and at the terminal time $T$, and in-between it is exponentially close to $0$. We have given explicit formulas, showing that the control is the sum of an exponentially decreasing term and of an exponentially increasing one.
If the time horizon is infinite (without final conditions), then only the first term remains, and the optimal control exponentially stabilizes the system, accordingly to the classical Riccati theory. The norms of the control action and of the optimal state decay then exponentially in time.
These results show that
as soon as the objective functional
of the optimal control problems for the considered system
contains a nontrivial tracking term,
the optimal solution has a special behavior referred to as the turnpike phenomenon.

Finally, we have shown that, if the final time $T$ is a positive even integer, then there exists a weight $\lambda$ such that the solution of the problem of exact controllability with minimal control norm coincides with the solution of the infinite horizon optimal control problem along the time interval $[0,2]$. This result justifies a receding horizon control strategy, where the first part of a finite horizon optimal control is used and then the procedure is updated in order to control the system over an infinite time horizon.

\medskip

As already said, the turnpike property has been much investigated in finite dimension (see \cite{trelat} and references therein for a general result).
In the infinite-dimensional setting, in
\cite{porretta} distributed control has been considered both for the heat equation
and the wave equation.

The turnpike phenomenon put in evidence in the present paper shows an interesting qualitative bifurcation of the HUM control as soon as the objective functional involves a tracking term. However, here, we have been able to show it by means of explicit computations.

Several open questions are in order.

First of all, it makes sense to consider an objective functional in which the tracking term is replaced with a discrepancy between the solution and a time-independent function, which is not necessarily a steady-state. According to the results of \cite{trelat}, we expect then that the turnpike property still holds true, and that, in large time, the optimal trajectory remains essentially close to the optimal steady state state, defined as the the closest steady-state to the objective.
However, in that case, we certainly do not have explicit formulas as derived in the present paper. Moreover, here we only considered a functional penalizing the normal derivative at $x=0$, and then we can only consider a time-independent function that is a steady-state, as said in Remark \ref{rem_steadystate}. But if instead, we were considering for instance the full norm in $H^1(0, 1)$, then we could consider in the objective functional a term of the form $\Vert y(t,\cdot)- a(\cdot)\Vert^2 $, where $a(\cdot)$ need not be a steady-state. Then, what can be expected is that, in large time, the optimal trajectory remains essentially close to the steady-state of the form $\sigma x$ that is the closest possible to the target $a(\cdot)$.

For more general multi-D wave equations, the situation is open.
Even if explicit computations can only be done in specific cases, we expect that the turnpike phenomenon is generic within the class of optimal control problems for controllable wave equations, and that HUM controls characterized by the adjoint system develop a quasi-periodic pattern, but when characterized by a more robust optimality cost, then, satisfy the turnpike property.

Another open issue is the investigation of semilinear wave equations (see \cite{coron}), for which steady-states may play an important role. Of course, in that case, we cannot expect that the turnpike property hold globally, but it should also hold as well at least in some neighborhood of an optimal steady-state (see discussions in \cite{trelat}).

\vspace*{0.5cm}

\bigskip
\noindent{\bf Acknowledgment.}
%This first author was
This work is  supported by DFG in the framework of the Collaborative Research Centre
CRC/Transregio 154, Mathematical Modelling, Simulation and Optimization Using the Example of Gas Networks, project C03 and project A03,
%
%The second and the third authors were partially supported by the Grant FA9550-14-1-0214 of the %EOARD-AFOSR.
%
%The third author was partially supported by the Grant MTM2011-29306-C02-00 of the MICINN (Spain), project PI2010-04 of the Basque Government, the ERC Advanced Grant FP7-246775 NUMERIWAVES, the ESF Research Networking Program OPTPDE.
%
%This work is supported
by the Advanced Grant FP7-246775 NUMERIWAVES of the European Research Council Executive Agency, FA9550-14-1-0214 of the EOARD-AFOSR, FA9550-15-1-0027 of AFOSR, the BERC 2014-2017 program of the Basque Government, the
MTM2011-29306-C02-00 and SEV-2013-0323 Grants of
the MINECO and a Humboldt Award at the University of Erlangen-Nuremberg.

%%%%%%%%%%%%%%%%%%%%%%%%%%%%%%%%%%%%%%%%%%%%%%%%%%%%%%%%%%%%%%%%%%%%%%%%%%%%

%\bibliographystyle{siam}

\end{document}